\documentclass[reqno,11pt]{amsart}

\usepackage[utf8]{inputenc}

\usepackage[T1]{fontenc}
\usepackage{lmodern}

\usepackage[a4paper]{geometry}

\DeclareUnicodeCharacter{00A0}{ }

\begin{document}
\title[Shape differentiation: non smooth case]
{Shape differentiation of \\ a steady-state reaction-diffusion problem \\ arising in Chemical Engineering: \\the case of non-smooth kinetic \\ with dead core}

\author[D. G\'omez-Castro\hfil EJDE-2017/??\hfilneg]
{David G\'omez-Castro}

\address{David G\'omez-Castro \newline
Instituto de Matem\'atica Interdisciplinar, Universidad Complutense de Madrid, Plaza de las Ciencias 3, 28040 Madrid, Spain}
\email{dgcastro@ucm.es}


\begin{abstract}
	In this paper we consider an extension of the results in shape differentiation of semilinear equations with smooth nonlinearity presented in D\'iaz, J.I., G\'omez-Castro, D.: An Application of Shape Differentiation to the Effectiveness of a Steady State Reaction-Diffusion Problem Arising in Chemical Engineering. Electron. J. Differ. Equations. 22, 31–-45 (2015) to the case in which the nonlinearities might be less smooth. Namely we will show that Gateaux shape derivatives exists when the nonlinearity is only Lipschitz continuous, and we will give a definition of the derivative when the nonlinearity has a blow up. In this direction, we will study the case of root-type nonlinearities. 
\end{abstract}

\maketitle
\numberwithin{equation}{section}
\newtheorem{theorem}{Theorem}[section]
\newtheorem{lemma}[theorem]{Lemma}
\newtheorem{remark}[theorem]{Remark}
\newtheorem{definition}[theorem]{Definition}
\newtheorem{corollary}[theorem]{Corollary}
\newtheorem{example}[theorem]{Example}
\newtheorem{proposition}[theorem]{Proposition}
\allowdisplaybreaks

\section{Introduction}

In this paper we consider the shape differentiation of a family of diffusion-reaction problems introduced by Aris in the context of optimization of chemical reactors depending on the spatial domain (see \cite{Aris+Strieder:1973}). It was later shown that the model can be rigorously deduced as a limit of different nonhomogeneous microscopic models (see  \cite{Conca+Diaz+Linan+Timofte:2004,diaz+gomez-castro+podolskii+shaposhnikova2017jmaa}). In particular we will be interested in the solutions of problem
\begin{equation}
\begin{cases}
-\Delta w+\beta (w)={f}, & \text{in }\Omega , \\ 
w=1, & \text{on }\partial \Omega ,%
\end{cases}
\label{eq:model problem w}
\end{equation}%
and their behaviour as we deform the domain $\Omega$.\\

It will be sometimes useful to consider the change in variable $u = 1 -w, g(u) = \beta(1) - \beta (1 - u) $ and $\widehat f = \beta (1) - f$, so that we have $u = 0$ on the boundary. After this change in variable we have that $u$ is the solution of
\begin{equation}
\begin{cases}
-\Delta u+g(u)=\widehat f, & \text{in }\Omega , \\ 
u=0, & \text{on }\partial \Omega .
\end{cases}%
\label{eq:model problem u}
\end{equation}%
These functions will be sometimes denoted $u_\Omega, w_\Omega$ when different domains are considered.\\

In \cite{Diaz+Gomez-Castro:2015shapediff} (see also \cite{Simon:1980differentiation,Henrot+Pierre:2005,Pironneau:1984}) the authors showed that, if $\beta \in W^{2,\infty} (\mathbb R)$ and $f \in L^2(\Omega)$ then the maps
\begin{eqnarray*}
	W^{1, \infty} (\mathbb R^n , \mathbb R^n ) & \to & H_0^1 (\Omega) \\
	\theta & \mapsto & u_{(I + \theta) \Omega} \circ ( I + \theta ) \\
	W^{1, \infty} (\mathbb R^n , \mathbb R^n ) & \to & L^2 (\mathbb R^n) \\
	\theta & \mapsto & u_{(I + \theta) \Omega}, 
\end{eqnarray*}
where the extension by $0$ is considered in $\mathbb R^n \setminus \Omega$, are Fréchet differentiable at $0$. Fixing $\theta \in W^{1,\infty} (\mathbb R^n, \mathbb R^n)$ it was shown in \cite{Diaz+Gomez-Castro:2015shapediff} that the directional derivative (the derivative of $u_\tau = u_{(I + \tau \theta) \Omega}$ with respect to $\tau$, $\frac{d u_\tau}{d\tau} = \frac{d u_\tau}{d\tau}|_{\tau = 0}$) is the solution of the following problem
\begin{equation}
\begin{cases}
-\Delta \frac{d u_\tau} {d\tau}+ g^{\prime }(u_\Omega)\frac{d u_\tau} {d\tau}=0, & 
\text{in }\Omega, \\ 
\frac{d u_\tau} {d\tau}=-\nabla u_\Omega \cdot \theta , & \text{on }%
\partial \Omega.%
\end{cases}
\label{eq:u' in terms of theta}
\end{equation}
Notice that, since $u = 1 -w$, we have that $\frac{d u_\tau} {d\tau} = -\frac{d w_\tau} {d\tau}$. Hence, taking into account that $g'(u) = -\beta' (w)$, we have that
\begin{equation}
\begin{cases}
-\Delta \frac{d w_\tau} {d\tau}+\beta^{\prime }(w_{\Omega})\frac{d w_\tau} {d\tau}=0, & 
\text{in }\Omega, \\ 
\frac{d w_\tau} {d\tau}=-\nabla w_\Omega \cdot \theta , & \text{on }%
\partial \Omega.%
\end{cases}
\label{eq:w' in terms of theta}
\end{equation}
The aim of this paper is to extend this kind of results to the case when $\beta \notin W^{2,\infty}$. First, we will show that, when $\beta \in W^{1,\infty}$ then the Gateaux shape derivative exists. 
However, if $\beta$ is not locally Lipschitz continuous,  the solution of \eqref{eq:model problem w} might develop a region of positive measure
\begin{equation}
\label{eq:dead core}
N_\Omega = \{ x \in \Omega: w_\Omega (x) = 0 \}.
\end{equation}
This region, known as \emph{dead core}\index{dead core}, was studied at length in \cite{Diaz:1985,bandle+sperb+stakgold1984}. It is a necessary condition for the existence of this region that $\beta' (w_\Omega) = +\infty$. Hence, equation \eqref{eq:w' in terms of theta} cannot be understood immediately in a standard way. In this setting, we will show that there exists a limit of the previous theory.

\section{Statement of results}

For the rest of the paper $\Omega \subset \mathbb R^n$ will be a fixed domain, of class $\mathcal C^2$, and $n \ge 2$.

\subsection{Existence and estimates of shape derivatives}

\subsubsection{Existence of Gateaux derivative when $\beta \in W^{1,\infty}$}

In \cite{Diaz+Gomez-Castro:2015shapediff} the authors prove the existence of a shape derivative in the Fréchet sense when $\beta \in W^{2,\infty}(\mathbb R)$. Nonetheless, as is it usually the case, the equation for the derivative is well defined in a straightforward way when $\beta \in W^{1, \infty}(\mathbb R)$. In fact, the following result shows that, if $\beta \in W^{1,\infty}(\mathbb R)$ rather than $W^{2,\infty}(\mathbb R)$, then the shape derivative exists only in the Gateaux sense, which is weaker than the Fréchet sense. 

\begin{theorem} \label{thm:existence of Gateuax derivate}
	Let $\theta \in W^{1, \infty} (\mathbb R^n, \mathbb R^n)$, $\beta \in W^{1,\infty} (\mathbb R)$ be nondecreasing such that $\beta (0) = 0$ and $f \in H^1 (\mathbb R^n)$.  Then, the applications 
	\begin{eqnarray*}
		\mathbb R &\to& L^2 (\Omega) \\
		\tau &\mapsto& u_{(I + \tau \theta) \Omega} \circ  (I + \tau \theta),
	\end{eqnarray*}
	and
	\begin{eqnarray*}
		\mathbb R &\to& L^2 (\mathbb R^n) \\
		\tau &\mapsto& u_{(I + \tau \theta) \Omega} 
	\end{eqnarray*}
	are differentiable at $0$. Furthermore, $\frac{d u_\tau}{d\tau}|_{\tau = 0}$ is the unique solution of \eqref{eq:u' in terms of theta}.
\end{theorem}

\begin{remark}
	In most case, the process of homogenization mentioned in the introduction gives an homogeneous equation \eqref{eq:model problem w} in which $\beta$ is the same as in the microscopic limit, and thus it is natural that $\beta$ be singular. However, it sometimes happens that the limit kinetic is different. 
	In the homogenization of problems with particles of critical size (see \cite{diaz+gomez-castro+podolskii+shaposhnikova2017generalmmg}) it turns out that the resulting kinetic in the macroscopic homogeneous equation \eqref{eq:model problem w} satisfies $\beta \in W^{1,\infty}$, even when the original kinetic of the microscopic problem was a general maximal monotone graph.
\end{remark}

\subsubsection{From $W^{2,\infty}$ to $W^{1,\infty} \cap \mathcal C^1$}

Let us show that the shape derivative is continuously dependent on the nonlinearity, and thus that we can make a smooth transition from the Fréchet scenario presented in \cite{Diaz+Gomez-Castro:2015shapediff} to our current case. For the rest of the paper we will use the notation:
\begin{equation}
v = \frac{d w_\tau}{d\tau}\Big|_{\tau = 0}
\end{equation}
\begin{lemma} \label{thm:approximation}
	Let $f \in L^2 (\mathbb R^n)$, $\beta \in W^{1,\infty} (\mathbb R)$ be nondecreasing functions such that $\beta(0) = 0$ and let $\beta_n \in W^{2,\infty} (\mathbb R)$ nondecreasing such that $\beta_n (0) = 0$. 
	Let $w_n$ be the unique solution of
	\begin{gather}
		\label{eq:model wn}
		\begin{cases}
			-\Delta w_n + \beta_n (w_n) = f & \Omega, \\
			w_n = 1 & \partial \Omega.
		\end{cases}
	\end{gather} 
	Then
	\begin{align}
		\| w_n - w \|_{H^1 (\Omega)} &\le C \| \beta_n - \beta \|_{L^\infty} \label{eq:H1 norm wm minus w}\\
		\| w_n - w\|_{H^2 (\Omega)} & \le C (1 + \| \beta' \|_{L^\infty} ) \| \beta_n - \beta \|_{L^\infty}. \label{eq:H2 norm wm minus w}
	\end{align}
	Furthermore, let $\beta \in C^1 (\mathbb R) \cap W^{1,\infty} (\mathbb R)$ and $v_n$ be the unique solution of
	\begin{gather}
		\label{eq:model vn}
		\begin{cases}
			-\Delta v_n + \beta_n' (w_n) v_n = 0 & \Omega, \\
			v_n + \nabla w_n \cdot \theta = 0 & \partial \Omega. 
		\end{cases}
	\end{gather}
	Then
	\begin{equation}
	v_n \rightharpoonup v \textrm{ in } H^1 (\Omega).
	\end{equation}
\end{lemma}

\begin{remark}
	In \eqref{eq:H1 norm wm minus w} the notation
	\begin{equation*}
		\| \beta_n - \beta \|_{L^\infty} = \sup_{x \in \mathbb R} |\beta_n (x) - \beta(x)|
	\end{equation*}
	does mean that either $\beta_n$ or $\beta$ are $L^\infty (\mathbb R)$ functions themselves, but rather that their difference is pointwise bounded, and, in fact, this bound is destined to go $0$ as $n\to +\infty$. We will use this notation throughout the paper.
\end{remark}

\subsubsection{Shape derivative with a dead core}

We can prove that the shape derivative in the smooth case has, under some assumptions, a natural limit when $\beta$ not smooth. \\

In some cases in the applications (see \cite{Diaz:1985}) we can take $\beta$ so that $\beta' (w_\Omega)$ has a blow up. It is common, specially in Chemical Engineering, that $\beta' (0) = +\infty$ and $N_\Omega$ exists (see \cite{Diaz:1985}). In this case $\beta' (w_\Omega) = +\infty$ in $N_\Omega$. Due to this fact, the natural behaviour of the weak solutions of \eqref{eq:w' in terms of theta} is $v = 0$ in $N_\Omega$. We have the following result

\begin{theorem} \label{thm:existence of derivative}
	Let $\beta$ be nondecreasing, $\beta(0) = 0$, $\beta' (0) = + \infty$,
	$$\beta \in \mathcal C(\mathbb R) \cap \mathcal C^1 (\mathbb R \setminus \{0\}),$$
	and assume that $|N_\Omega|> 0$, $\theta \in W^{1,\infty} (\mathbb R^n, \mathbb R^n)$ and
	$
	0 \le f \le \beta(1).
	$
	Then, there exists $v$ a solution of
	\begin{equation} \label{eq:model v}
	\begin{cases}
	-\Delta v + \beta'(w_\Omega) v = 0 & \Omega \setminus N_\Omega, \\
	v = 0 & \partial N_\Omega, \\
	v = - \nabla w_\Omega \cdot \theta & \partial \Omega,
	\end{cases}
	\end{equation}
	in the sense that $v \in H^1 (\Omega)$, $v = 0$ in $N_\Omega$, $v = - \nabla w_\Omega \cdot \theta$ in $L^2 (\partial \Omega)$, $\beta' (w_\Omega) v^2 \in L^1 (\Omega)$ and 
	\begin{equation}
	\int _{\Omega \setminus N_\Omega} \nabla v \nabla \varphi + \int _{\Omega\setminus N_\Omega} \beta'(w) v \varphi = 0 
	\end{equation}
	for every $\varphi \in W^{1, \infty}_c (\Omega \setminus N_\Omega) $.
	Furthermore, for $m \in \mathbb N$, consider $\beta_m$ defined by
	$$\beta_m'(s) = \min \{  
	m , \beta' (s) \}, 
	\qquad \beta_m(0) = \beta(0) = 0, $$
	and let $w_m, v_m$ be the unique solutions of \eqref{eq:model wn} and \eqref{eq:model vn}.
	Then,
	\begin{equation}
	v_m \rightharpoonup v, \quad \textrm{in }H^1 (\Omega),
	\end{equation}
	where $v$ is a solution of \eqref{eq:model v}.
\end{theorem}

The uniqueness of solutions of \eqref{eq:model v} when $\beta'(w_\Omega)$ blows up is by no means trivial. Problem \eqref{eq:model v} can be written in the following way:
\begin{equation}
\label{eq:linear problem with potential}
-\Delta v + V v  = f 
\end{equation} 
where $V = \beta' (w_\Omega)$ may blow up as a power of the distance to a piece of the boundary. This kind of problems are common in Quantum Physics, although their mathematical treatment is not always rigorous (cf.  \cite{Diaz2015ambiguousSchrodinger,Diaz2017ambiguousSchrodinger}). \\

In the next section we will show estimates on $\beta' (w_\Omega)$. Let us state here some uniqueness results depending on the different blow-up rates.\\

When the blows is subquadratic (i.e. not \emph{too} rapid), by applying Hardy's inequality and the Lax-Migram theorem, we have the following result (see \cite{Diaz2015ambiguousSchrodinger,Diaz2017ambiguousSchrodinger}).
\begin{corollary}
	Let $N_\Omega$ have positive measure and $\beta'(u(x)) \le C d(x,  N_\Omega)^{-2}$ for a.e. $x \in \Omega \setminus N_\Omega$. Then the solution $v$ is unique.
\end{corollary}

The study of solutions of problem \eqref{eq:linear problem with potential} in $\Omega$ when $V \in L^1_{loc} (\Omega)$ by many authors (see \cite{Diaz+Rakotoson:2010,Diaz+Gomez-Castro:2015veryweak} and the references therein). Existence and uniqueness of this problem in the case $V(x) \ge C d (x, \partial \Omega)^{-r}$ with $r > 2$ was proved in \cite{Diaz+Gomez-Castro:2015veryweak}. Applying these techniques one can show that
\begin{corollary}
	Let $N_\Omega$ have positive measure and $\beta'(u(x)) \ge C d(x,  N_\Omega)^{-r}, r > 2$ for a.e. $x \in \Omega\setminus N_\Omega$. Then the solution $v$ is unique.
\end{corollary}
Similar techniques can be applied to the case $\beta'(u(x)) \ge C d(x,  N_\Omega)^{-2}$. This will be the subject of a further paper.

\subsection{Estimates of $w_\Omega$ close to $N_\Omega$}

Let us study the solution $w_\Omega$ on the proximity of the dead core and the blow up behaviour of $\beta' (w_\Omega)$. First, we present a known example
\begin{example}
	Explicit radial solutions with dead core are known when $\beta (w) = |w|^{q-1}w$ ($0 < q	< 1$), $\Omega$ is a ball of large enough radius and $f$ is radially symmetric. In this case it is known that $N_\Omega$ exists, has positive measure and
	$$
	\frac{1}{C} d(x,  N_\Omega)^{-2}\le \beta' (w_\Omega) \le C d(x,  N_\Omega)^{-2}.
	$$
	For the details see \cite{Diaz:1985}.
\end{example}

In fact, we present here a more general result to study the behaviour in the proximity of the dead core, based on estimates from \cite{Diaz:1985}.

\begin{proposition} \label{prop:general nonlinearity}
	Let $f = 0$, $\beta$ be continuous, monotone increasing such that $\beta(0) = 0$, $w$ be a solution of \eqref{eq:model problem w} that develops a dead core $N_\Omega$ of positive measure and $\partial N_\Omega \in \mathcal C^1$. Assume that
	\begin{equation}
	G(t) =  \sqrt 2 \left(  \int_0^t \beta (\tau) d \tau + \alpha t  \right)^{\frac 1 2}, \quad \textrm{ where } \alpha =  \max \left\{  0 , \min_{x\in\partial \Omega} H(x) \frac{\partial w}{\partial n}(x)  \right\} ,
	\end{equation}
	is such that $\frac 1 G \in L^1(\mathbb R)$. Then 
	\begin{equation}
	w_\Omega(x) \le \Psi^{-1} (d(x,{N_{\Omega}})), \quad \textrm{ where } \Psi (s) =  \int_{0}^s \frac{dt} { G(t) },
	\end{equation}
	in a neighbournood of $N _\Omega$.
\end{proposition}

\begin{example}[Root type reactions]
	Let $f = 0$, $\beta(s) = \lambda |s|^{q-1} s$ with $0 < q < 1$ and $\Omega$ be convex such that $N_\Omega$ exists and $\partial N_\Omega \in \mathcal C^1$. Then
	\begin{equation}
	w_\Omega(x) \le C d(x, {N_\Omega} )^{\frac 2 {1-q}}.
	\end{equation}	
	Furthermore
	\begin{equation}
	\beta'(w_\Omega (x)) \ge Cd(x, {N_\Omega}) ^{ -2 }.
	\end{equation}
\end{example}


\section{Proof of Theorem \ref{thm:existence of Gateuax derivate}}
For the rest of the paper let us note 
\begin{equation} 
u_\tau = u_{(I + \tau \theta) \Omega}.
\end{equation}
Notice that $u_0 = u_\Omega$.

Let us define $U_\tau = u_{(I + \tau \theta)\Omega} \circ (I + \tau \theta) \in H_0^1 (\Omega)$. Again $U_0 = u_0 = u_\Omega$. We have that
\begin{equation}
\int _\Omega A_\tau \nabla U_\tau \nabla \varphi  + \int _ \Omega g (U_\tau) \varphi J_\tau = \int _ \Omega f_\tau  \varphi J_\tau ,
\end{equation}
where $J_\tau$ is the Jacobian of the transformation. $f_\tau = f \circ (I + \tau \theta)$ and $A_\tau$ is the corresponding diffusion matrix (see \cite{Diaz+Gomez-Castro:2015shapediff} for the explicit expression). Fortunately, $J_\tau \ge 0$ and, for $\tau$ small, we have that $\xi \cdot A_\tau \xi  \ge A_0 |\xi|^2$ for some $A_0 > 0$ constant.
Considering the difference of the weak formulations of $U_\tau$ and $U_0 = u_\Omega$ we have that
\begin{align*}
	\int_{\Omega} A_\tau \nabla (U_\tau - u_0) \nabla \varphi + \int _ \Omega (g(U_\tau) - g(u_0) ) J_\tau \varphi & = \int _ \Omega (f_\tau J_\tau - f) \varphi + 
	\\
	&\qquad + \int _ \Omega (I - A_\tau) \nabla u_0 \nabla \varphi  
	\\
	& \qquad
	+ \int _ \Omega (J_\tau -1 ) g(u_0) \varphi. 
\end{align*}
Hence, due to the monotonicity of $g$, we have that
\begin{equation*}
	\left \| \nabla \left( \frac{U_\tau - u}{\tau } \right) \right \|_{L^2  } \le C \left ( \left\| \frac{f_\tau - f}{\tau}  \right\|_{L^2} + \left\| \frac{ A_\tau - I}{\tau }   \right\| _{L^\infty} \| \nabla u_0 \|_{L^2} + \left\|  \frac{ J _ \tau - 1}{\tau}  \right\|_{L^\infty} \|g(u_0)\|_{L^2} \right )
\end{equation*}
Since $f_\tau, A_\tau$ and $J_\tau$ are differentiable at $0$,  there is weak $H_0^1 (\Omega)$ limit. Hence, the limit is strong in $L^2 (\Omega)$. Therefore, the function
\begin{equation}
u_\tau = U_\tau \circ (I + \tau \theta )^{-1}
\end{equation}
is differentiable with respect to $\tau \in \mathbb R$ with images in $L^2 (\Omega)$ at $\tau = 0$. Besides, 
\begin{equation}
H_0^1 (\Omega) \ni \frac{d U_\tau}{d \tau} \Big|_{\tau = 0}  = \frac{du_\tau}{d\tau} \Big|_{\tau = 0} + \nabla u_0 \cdot \theta.
\end{equation}
To characterize the derivative, we differenciate on the
variational formulation
\begin{equation*}
	\int_{\mathbb{R}^n} f \varphi   = \int_{\mathbb{R}^n} \left( -u_\tau \Delta \varphi +
	g(u_\tau) \varphi \right)  \qquad \forall \varphi \in \mathcal C_c ^\infty (\Omega).
\end{equation*}
Considering the difference of the equations for $u_\tau$ and $u_0$ and diving by $\tau$
\begin{align}
	0 &= \int _{\mathbb R^n } \left( -\frac{u_\tau - u_0}{\tau } \Delta \varphi + \frac{g(u_\tau) - g(u_0)}{\tau} \varphi \right) \\
	&= \int _{\mathbb R^n } \frac{u_\tau - u_0}{\tau } \left(  - \Delta \varphi + \frac{g(u_\tau) - g(u_0)}{u_\tau - u_0} \varphi  \right).
\end{align}
Notice that
\begin{equation*}
	\left| \frac{g(u_\tau) - g(u_0)}{u_\tau - u_0} \right| \le \| g ' \|_{L^\infty}.
\end{equation*}
Therefore, up to a subsequence, $\frac{g(u_\tau) - g(u_0)}{u_\tau - u_0}$ converges weakly in $L^2 (\Omega)$. On the other hand since $u_\tau \to u_0$ pointwise, again up to a subsequence, so 
\begin{equation}
\frac{g(u_\tau) - g(u_0)}{u_\tau - u_0} \to g' (u_0) \quad \textrm{ a.e. in } \Omega.
\end{equation}
Via a Césaro mean argument we have that the weak $L^2$ limit and pointwise limit coincide. Hence, passing to the limit in $L^2 (\Omega)$
\begin{equation}
0 = \int_{\Omega} \frac{du_\tau}{d\tau }\Big|_{\tau = 0}\left( -  \Delta \varphi +  g^{\prime }(u_0)
\varphi \right), \qquad \varphi \in \mathcal C_c ^\infty (\Omega).
\end{equation}
Therefore $\frac{du_\tau}{d\tau}$ is the unique solution of \eqref{eq:u' in terms of theta}. \qed

\section{Proof of Lemma \ref{thm:approximation}}
By considering the difference of the weak formulations we have that
\begin{equation*}
	\int_ \Omega \nabla (w_m - w) \nabla \varphi + \int_{\Omega } (\beta_m (w_m) - \beta_m (w)) \varphi = \int_{ \Omega } (\beta(w) - \beta_m (w) ) \varphi .
\end{equation*}
Taking $\varphi = w_m - w$, and using the monotonicity of $\beta_m$ we have that
\begin{equation*}
	\| \nabla (w_m - w) \|_{L^2} ^2  \le \| \beta_m - \beta \|_{L^\infty} \| w_m - w \|_{L^1 (\Omega)}.
\end{equation*}
Using Poincaré inequality and the embedding $L^1 \hookrightarrow L^2$ we have that
\begin{equation*}
	\| w_m - w \|_{L^2} \le C \| \beta_m - \beta \|_{L^\infty}.
\end{equation*}
By considering the equation
\begin{align*}
	\| \Delta (w_m - w) \|_{L^2} &= \| \beta(w) - \beta_m (w_m) \|_{L^2} 
	\\
	&\le \| \beta(w) - \beta(w_m) \|_{L^2} + \| \beta(w_m) - \beta_m (w_m) \|_{L^2} \\
	&\le \| \beta' \|_{L^\infty } \| w_m - w \|_{L^2} + \| \beta_m - \beta \|_{L^\infty}.
\end{align*}
Hence, to deduce \eqref{eq:H2 norm wm minus w} we apply that
\begin{align*}
	\| w_m - w \|_{H^2 } &\le C ( \| \Delta( w_m - w ) \|_{L^2 } + \|w_m - w \|_{L^2}).
\end{align*}
Considering the difference of the weak formulations of the problems for $v_m$ and $v$ we have that
\begin{align}
	\int _ \Omega \nabla (v_m - v) \nabla \varphi &= \int _ \Omega (\beta'(w) v - \beta_m'(w_m) v_m) \varphi 
	\nonumber 
	\\ 
	& = \int _\Omega (\beta' (w) - \beta_m'(w_m)) v_m \varphi + \int _ \Omega \beta' (w)(v - v_m) \varphi 
	\nonumber 
	\\
	& = \int_ \Omega (\beta' (w) - \beta' (w_m)) v_m \varphi + \int _ \Omega (\beta' (w_m) - \beta_m' (w_m)) v_m \varphi \nonumber
	\\
	& \quad + \int _ \Omega \beta' (w) (v - v_m) \varphi  \label{eq:weak formulation difference vs} 
\end{align}
for all $\varphi \in H_0^1 (\Omega)$. Considering the test function $\varphi = v_m - v + \nabla (w_m - w) \cdot \theta \in H_0^1 (\Omega)$ we have, applying \eqref{eq:H2 norm wm minus w}
\begin{align*}
	\int _ \Omega |\nabla (v_m - v)|^2 &\le C (1 + \| w_m- w\|_{H^2}) 
	\\
	&\qquad   \times \Big( (1 + \| \beta' (w)\|_{L^\infty})\| w_m- w\|_{H^2} \\
	&\qquad \qquad  + \|v_m\|_{L^2} (\| \beta_m' + \beta' \|_{L^\infty} + \| \beta' (w_m) - \beta'(w)\|_{L^\infty}   )\Big )
	.
\end{align*}
We cannot guaranty that $\|\beta' (w_m) - \beta' (w)\|_\infty $ goes to zero. However it is, indeed, bounded by $2 \| \beta' \|_{L^\infty}$.
On the other hand, taking into account the boundary condition
\begin{equation}
\| v_m - v \|_{L^2 (\partial \Omega)} \le C \| \nabla (w_m - w)\|_{L^2 (\partial \Omega)} \le C \| w_m - w \|_{H^2 (\Omega)} \le C \| \beta_m - \beta \|_{L^2} \to 0. \label{eq:convergence of boundary condition difference vs}
\end{equation}
Hence, there is a weak limit $\widehat v \in H^1(\Omega)$
\begin{equation}
v_m - v \rightharpoonup \widehat v \textrm{ in } H^1 (\Omega).
\end{equation}
Due to \eqref{eq:convergence of boundary condition difference vs} we have that $\widehat v \in H_0^1 (\Omega)$. Taking into account \eqref{eq:weak formulation difference vs} and the fact that $\beta'(w_m) \to \beta' (w)$ a.e. in $\Omega$, have that
\begin{equation}
\int_\Omega \nabla \widehat v \nabla \varphi + \int_ \Omega \beta' (w) \widehat v \varphi = 0 \qquad \forall \varphi \in H_0^1 (\Omega). 
\end{equation}
Taking $\varphi = \widehat v \in H_0^1 (\Omega)$ as a test function we deduce that $\widehat v = 0$.
\qed

\section{Proof of Theorem \ref{thm:existence of derivative}}
We start by pointing out that, due to the condition on $f$ we have that $0 \le  w_m \le 1$.
Since $\beta_m \nearrow \beta$ in $[0,1]$ we have $w_m$ is pointwise decreasing (see \cite{Evans1998}). 
Hence, there exists a pointwise limit $w$ such that $w_m \searrow w$ a.e. in $\Omega$. In particular $0 \le w \le 1$. 
Due to the Dominated Convergence Theorem  we have that
\begin{equation}
w_m \to w \textrm{ in } L^p (\Omega) \quad \forall 1 \le p < +\infty.
\end{equation} 
Let $U \subset \Omega$ be an open neighbourhood of $\partial \Omega$ such that $\overline U \cap N_\Omega = \emptyset$ and $\partial U \in \mathcal C^2$. Then
\begin{equation}
\underline w_U = \inf _{U} w > 0.
\end{equation}
We have that $w_m \ge w \ge  \underline w_U$. We have that $\beta  \in \mathcal C^1([\underline w_U, 1])$ and, hence, $\beta_m \to \beta$ in $\mathcal C^1([\underline w_U, 1])$. Therefore 
\begin{equation}
\beta_m(w_m) \to \beta(w) \textrm{ in } L^p (\Omega \setminus \overline U) \quad \forall 1 \le p < +\infty, \\
\end{equation}
Since $\| w_m\|_{H^1} \le C(1 + \| \beta_m (w_m)\|_{L^2} + \| f\|_{L^2} )$ we have that $w_m \rightharpoonup w$ in $H^1 (\Omega)$, and thus that $w$ is the unique solution of \eqref{eq:model problem w}.  
Applying this
\begin{equation}
\Delta w_m = \beta_m (w_m) - f \to \beta (w) - f = \Delta w  \textrm{ in } L^p (\Omega \setminus \overline U).
\end{equation}
Thus 
\begin{equation}
\| w_m - w\| _{H^2 (\Omega \setminus \overline U)} \le C (\| \Delta (w_m - w) \|_{L^2 (\Omega \setminus \overline U)} + \|w_m - w\|_{L^2(\Omega \setminus \overline U)}) \to 0.
\end{equation}
Hence
$$ w_m \to w \textrm{ in } H^2(\Omega \setminus \overline U). $$		
In particular
$$		\nabla w_m \to \nabla w \textrm{ in } H^{\frac 1 2}(\partial \Omega)^n. $$
Since $\beta_m' \in L^\infty (\mathbb R)$ we take the ``shape derivative'' $v_m$ solution of \eqref{eq:model vn}, which is well defined. Let us find their limit.
\\
Let us show we show that 
\begin{equation}
\beta_m' (w_m ) \to \beta' (w) \textrm{ a.e. in } \Omega.
\label{eq:pointwise convergence potential}
\end{equation}
First, let $x \notin N_\Omega$. Then $\beta$ is $C^1$ in $w(x)$. Therefore $\beta' (w_m(x)) \to \beta' (w(x))$. Hence, the sequence $\beta' (w_m (x))$ is bounded, so $ \beta' (w_m (x)) \le m_0$ for some $m_0$ large. Thus
$
\beta_m' (w_m (x)) = \beta' (w_m (x) ) 
$	
for $m \ge m_0$. Hence the convergence is proved for $x \notin N_\Omega$. Let $x \in N_\Omega$. Then $\beta' (w(x)) = +\infty$. Since $w_m (x) \to w (x)$ then $\beta' (w_m (x)) \to +\infty$. In that case, we have that 
\begin{equation*}
	\beta_m' (w_m (x)) = \beta (w_m (x)) \wedge m \to + \infty = \beta (w(x)).
\end{equation*}
This completes the proof of \eqref{eq:pointwise convergence potential}.
\\
Let us show that sequence $(v_m)$ is bounded in $H^1 (\Omega )$. There exist two open sets $U_0, U_1 \subset \Omega$ such that $\partial \Omega \subset U_1, N_\Omega \subset U_0$, $U_0 \cap U_1 = \emptyset$. There also exists a smooth transition function $\Psi$ such that $\Psi = 0$ in $U_0$ and $\Psi = 1$ in $U_1$. Let us define $g_m = \Psi \nabla w_m \cdot \theta \in H^1(\Omega)$. Then $\varphi = v_m + g_m \in H_0^1(\Omega)$ and it can be used as a test function in the weak formulation. Hence
\begin{align*}
	\int _ \Omega \nabla v_m \nabla (v_m + g_m) + \int _ \Omega \beta_m' (w_m) v_m (v_m + g_m) = 0 . 
\end{align*}
Therefore, through standard arguments
\begin{align*}
	\int_ \Omega |\nabla v_m|^2 + \int_ \Omega \beta_m'(w_m) v_m^2 
	&= - \int _ \Omega \nabla v_m \nabla g_m - \int _\Omega \beta'_m(w_m) v_m g_m \\ 
	&\le  \left( \int _ \Omega |\nabla v_m|^2 \right)^{\frac 1 2}  \left( \int _ \Omega |\nabla g_m|^2 \right)^{\frac 1 2} 
	\\
	&\qquad + \left( \int _\Omega \beta'_m(w_m) v_m^2 \right)^{\frac 1 2} \left( \int _\Omega \beta'_m(w_m) g_m^2 \right)^{\frac 1 2}
	\\
	& \le \frac{1}{2} \left( \int_ \Omega |\nabla v_m|^2 + \int_ \Omega \beta_m'(w_m) v_m^2  \right) 
	\\
	&\qquad + C \left(  \int_ \Omega |\nabla g_m|^2 + \int_ \Omega \beta_m'(w_m) g_m^2  \right).
\end{align*}
Since $\beta'_m(w_m)$ is uniformly bounded in $L^\infty (\Omega \setminus \overline {U_0})$ we have that the sequence is bounded:
\begin{align*}
	\left( \int_ \Omega |\nabla v_m|^2 + \int_ \Omega \beta_m'(w_m) v_m^2  \right) & \le C \left(  \int_ \Omega |\nabla g_m|^2 + \int_ \Omega \beta_m'(w_m) g_m^2  \right) \le C.
\end{align*}
In particular, there exists $v \in H^1 (\Omega)$ such that, up to a subsequence,
\begin{equation*}
	v_m \rightharpoonup v \textrm{ in } H^1 (\Omega).
\end{equation*}
Also, due to Fatou's lemma
\begin{equation}
\int_ \Omega \beta'(w) v^2 \le C.
\end{equation}
Since $\beta' (w) = + \infty$ in $N_\Omega$ we have that $v = 0$ a.e. in $N_\Omega$. 
For $\varphi \in W_c^{1,\infty} (\Omega \setminus N_\Omega)$ we have that
\begin{equation}
\int_{\Omega \setminus N_\Omega} \nabla v_m \nabla \varphi +\int_{\Omega \setminus N_\Omega} \beta'_m(w_m) v_m \varphi = 0.
\end{equation}
Let us consider the compact subset $K = \textrm{supp} \varphi \subset \Omega \setminus N_\Omega$.
Let us show that $\beta'(w_m) \to \beta'(w)$ in $L^2 (K)$. We have  $0 < \underline w_K \le w \le w_m$ in $K$. Due to the Dominated Convergence Theorem we have that $\beta_m'(w_m) \to \beta' (w)$ strongly in $L^p (K)$ for $1 \le p < + \infty$. 
\\
Hence, by passing to the limit we deduce that
\begin{equation}
\int_{\Omega \setminus N_\Omega} \nabla v \nabla \varphi +\int_{\Omega \setminus N_\Omega} \beta'(w) v \varphi = 0.
\end{equation}
This completes the proof.
\qed

\section{Proof of Proposition \ref{prop:general nonlinearity}}

Let us consider $x_0 \in \partial N_\Omega$ and
\begin{equation}
W(t) = w_\Omega(x_0 + t n (x_0))
\end{equation}
where $n(x_0)$ represents the normal vector to $ \partial N_\Omega$ at $x_0$. Due to Theorem 1.24 in \cite{Diaz:1985}, we have that
\begin{equation}
\frac 1 2 |\nabla w_\Omega (x)|^{2} \le \int_0^{w_\Omega(x)} \beta(s)ds + \alpha w_\Omega(x)
\end{equation}
for all $x \in \overline \Omega$. Hence
\begin{align*}
	\frac{d W}{dt} &\le \left| \frac{d W}{dt} \right| = | \nabla w_\Omega (x_0 + t n(x_0)) \cdot n(x_0)| 
	\\
	&\le |\nabla w_\Omega (x_0 + t n(x_0)) | \le G(w_\Omega (x_0 + t n(x_0))) 
	\\
	&= G (W(t)).
\end{align*}
Thus, $W$ is a solution of the following Ordinary Differential Inequality
\begin{equation}
\begin{cases}
\frac{d W}{dt} (t) \le G(W(t)), \\
W(0) = 0.
\end{cases}
\end{equation}
Let us consider $W_\varepsilon$ the solution of 
\begin{equation}
\begin{cases}
\frac{d W_\varepsilon}{dt} (t) = G(W_\varepsilon(t)), \\
v_\varepsilon(0) = \varepsilon.
\end{cases}
\end{equation}
This problem has a unique smooth solution, since $G \in \mathcal C^1(\mathbb R\setminus \{0\})\cap \mathcal C(\mathbb R)$ is strictly increasing and $G(0) = 0$. In fact, solving this simply separable O.D.E., we obtain that
\begin{equation} \label{eq:defn W eps}
W_\varepsilon (t) = \Psi^{-1} ( t + \Psi (\varepsilon)).
\end{equation}
Due to the monotonicity of $G$ we have that
\begin{equation}
W(t) \le W_\varepsilon (t) \quad \forall t \ge 0.
\end{equation}
Passing to the limit as $\varepsilon \to 0$ in \eqref{eq:defn W eps} we have that
\begin{equation}
	W(t) \le  \Psi^{-1} ( t ).
\end{equation}
Hence, since we can parametrize a neighbourhood of $\partial N_\Omega$ by $(x, t) \in \partial N_\Omega \times (-\lambda_0, \lambda_0) \mapsto x + t n (x) $, we deduce that
\begin{equation}
w(x) \le \Psi ^{-1} ( d (x, N _\Omega)	)
\end{equation}
at least in a neighbournood of $\partial N_\Omega$. This proves the result.	\qed

\section*{Acknowledgments}
The author is thankful to Professor Jes\'us Ildefonso D\'iaz for fruitful discussions in the preparation of this paper and his continued support. The research of D. G\'omez-Castro was supported by the Spanish goverment through an FPU fellowship (ref. FPU14/03702) and by the project ref. MTM2014-57113-P of the DGISPI.


\end{document}